\documentclass[12pt,leqno]{amsart}
\usepackage{amssymb}
\usepackage{amscd}
\usepackage[all]{xypic}

\newcounter{alphalistcntr}

\theoremstyle{plain}
\newtheorem{theorem}{Theorem}[section]
\newtheorem{lemma}[theorem]{Lemma}
\newtheorem{corollary}[theorem]{Corollary}
\newtheorem{prop}[theorem]{Proposition}

\theoremstyle{remark}

\newtheorem{remark}[theorem]{Remark}

\newtheorem*{note*}{Note}
\newtheorem*{remark*}{Remark}
\newtheorem*{example*}{Example}

\theoremstyle{definition}
\newtheorem*{definition*}{Definition}
\newtheorem{definition}[theorem]{Definition}

\newcommand{\Z}{\mathbb{Z}}

\newcommand{\Q}{\mathbb{Q}}

\newcommand{\N}{\mathbb{N}}
\newcommand{\Primes}{\mathbb{P}}

\newcommand{\Gal}{\mathrm{Gal}}
\newcommand{\tensor}{\otimes}
\newcommand{\Span}{\mathrm{Span}}
\newcommand{\rank}{\mathrm{rank}}

\title[On the trace map between abelian number fields]{On the trace map between absolutely abelian number fields of equal conductor}
\author{Henri Johnston}
\thanks{Partially supported by the States of Jersey Education, Sport and Culture Committee through the Jersey Scholarship.}
\address{Department of Mathematics \\ Malott Hall \\ Cornell University \\
Ithaca, NY 14853-4201 \\ USA}
\email{henri@math.cornell.edu}
\urladdr{http://www.math.cornell.edu/$\sim$henri}

\begin{document}

\maketitle

\section{Introduction}

Let $L/K$ be an extension of absolutely abelian number fields of equal conductor, $n$.
If $T_{L/K}: L \rightarrow K$ denotes the trace map, then $T_{L/K}(\mathcal{O}_{L})$ is an ideal in $ \mathcal{O}_{K}$. Let $I(L/K)$ denote the norm of $T_{L/K}(\mathcal{O}_{L})$ over $\Q$, i.e. 
$[\mathcal{O}_{K} : T_{L/K}(\mathcal{O}_{L})]$. Sharpening the main result of Girstmair in \cite{girstmair}, we determine $I(L/K)$ exactly for any such $L/K$: if $e=v_{2}(n)$ and $m=n/2^{e}$, then 
$$
I(L/K) = \left\{ 
\begin{array}{ll} 
2^{[K \cap \Q^{(m)}:\Q]} =  2^{( [K:\Q] / 2^{(e-2)} )} & \textrm{ if $L/K$ is wildly ramified,} \\
1 & \textrm{ otherwise.}
\end{array} \right.
$$
After first determining criteria for wild ramification of $L/K$ (which can only happen at primes above $2$), 
the above result 
is obtained for $n=2^{e}$ ($e \geq 3$) by computing $T_{L/K}(\mathcal{O}_{L})$ explicitly, and is then extended to the general case. This approach does not rely on Leopoldt's Theorem, in contrast to the techniques used 
in \cite{girstmair}. \par 
The explicit nature of the calculations used to compute $I(L/K)$ leads to the definition of an ``adjusted trace map'' $\hat{T}_{\Q^{(n)}/K}$ with the property that $\hat{T}_{\Q^{(n)}/K}(\mathcal{O}^{(n)})=\mathcal{O}_{K}$ (here $\Q^{(n)}$ denotes the $n^{\mathrm{th}}$ cyclotomic field and $\mathcal{O}^{(n)}$ its ring of integers). Using this map, we restate Leopoldt's Theorem and show that its proof can be reduced to the (easier) cyclotomic case.

\section{Dirichlet Characters}

We first recall some basic facts about Dirichlet characters. For more details, see Chapter 3 of \cite{wash} and Section 2 of \cite{lettl}.

\begin{definition}
For $n \in \N$, let $\zeta_{n}$ be a primitive $n^{\textrm{th}}$ root of unity and $\Q^{(n)}=\Q(\zeta_{n})$ the $n^{\textrm{th}}$ cyclotomic field. Let $\mathcal{O}^{(n)}=\mathcal{O}_{\Q^{(n)}}=\Z[\zeta_{n}]$ denote the ring of integers of $\Q^{(n)}$, and 
$X^{(n)}$ denote the group of Dirichlet characters of conductor dividing $n$. \par
Let $\Primes$ denote the set of rational primes. 
Define $p^*=4$ if $p=2$, and $p^*=p$ if $p \in \Primes, p \neq 2$.
\end{definition}

\begin{prop}\label{decomp}
Let $p \in \Primes$ and $e \in \N$, with $e \geq 2$ if $p=2$. Then there exists a natural decomposition
$$(\Z / p^e\Z)^{\times} = (\Z / p^{*}\Z)^{\times} \times (1+p^{*} \Z) / (1+p^e \Z) $$
where both factors are considered as subgroups of $(\Z / p^e\Z)^{\times}$.
Note that we take $(\Z / 4\Z)^{\times} = \{ \pm  1 \}$.
\end{prop}

\begin{proof}
Straightforward.
\end{proof}

\begin{definition}
Let $p \in \Primes$ and $e \in \N$ with $e \geq 2$ if $p=2$. Then dualizing the decomposition of Proposition \ref{decomp} yields the decomposition 
$$X^{(p^e)} = \langle \omega_p \rangle \times \langle \psi_{p^e} \rangle$$ with $\langle \omega_p \rangle  = X^{(p^*)}$ and $ \langle \psi_{p^e} \rangle$ the group of Dirichlet characters whose conductors divide $p^e$ and which are trivial on the factor $(\Z / p^*\Z)^{\times}$.
\end{definition}

\begin{theorem}
Let $n \in \N$. There is an order preserving one-to-one correspondence between subgroups of
$X^{(n)}$ and subfields of $\Q^{(n)}$. Let $X_i$ be the subgroup corresponding to the subfield $K_i$. 
Then $|X_i|=[K_i: \Q]$ and the compositum $K_1K_2$ corresponds to $\langle X_1, X_2 \rangle$.
\end{theorem}

\begin{proof}
See Chapter 3 of \cite{wash}.
\end{proof}

\begin{definition}
Let $p \in \Primes$, $X \subseteq X^{(n)}$ and $e=v_{p}(n)$. Then $X_{p}$ denotes the image of $X$ under the natural projection $\pi_{p}:X^{(n)} \rightarrow X^{(p^{e})}$.
\end{definition}

\begin{theorem}\label{ramindex}
Let $X$ be a group of Dirichlet characters and let $K$ be the associated abelian number field. 
Then $p \in \Primes$ has ramification index  $|X_{p}|$ in $K$.
\end{theorem}

\begin{proof}
This is Theorem 3.5 of \cite{wash}.
\end{proof}

\begin{remark}\label{wildtameparts}
When $p$ is odd, $\langle \omega_p \rangle$ and  $\langle \psi_{p^e} \rangle$ have orders $p-1$
and $p^{(e-1)}$ respectively. So by considering the decomposition $X^{(p^e)} = \langle \omega_p \rangle \times \langle \psi_{p^e} \rangle$, the field corresponding to $\langle \omega_p \rangle$  can be thought of as the ``tame part'' of $\Q^{(p^{e})}$, and that corresponding to $\langle \psi_{p^e} \rangle$ as the ``wild part''. \par
When $p=2$,
$\langle \omega_{2} \rangle$ and $\langle \psi_{2^{e}} \rangle$ have orders $2$ and $2^{(e-2)}$ respectively, and therefore both correspond to wildly ramified extensions of $\Q$ (namely $\Q(i)$ and the maximal totally real subfield
$\Q(\zeta_{2^{e}}+\zeta_{2^{e}}^{-1})$, respectively). In other words, $\Q^{(2^{e})}$ has no ``tame part''.
\end{remark}

\begin{prop}\label{surjective}
Let $K / \Q$ be an abelian extension of conductor $n=p_1^{e_1} \cdots p_t^{e_t}$ where $p_{1}=2$,
and let $X \subseteq X^{(n)}$ be its associated group of Dirichlet characters.
\begin{enumerate}
\item The natural projection
$ \pi_{\psi}: X \longrightarrow \prod_{i=1}^t \langle \psi_{p_i^{e_i}} \rangle$
is surjective.
\item Let $e=e_{1}=v_{2}(n)$. Then $X_{2}$ is either $X^{(2^{e})}=\langle \omega_{2} \rangle \times \langle \psi_{2^{e}} \rangle$, $\langle \psi_{2^{e}} \rangle$, or $\langle \omega_{2} \psi_{2^{e}} \rangle$. Note that $\psi_{2^{e}}$ is trivial if $e \leq 2$.
\item $\langle X, \prod_{i=2}^{t} \langle \omega_{p_{i}} \rangle \rangle = X_{2} \times 
\prod_{i=2}^{t} \langle \omega_{p_i} \rangle \times \prod_{j=2}^{t} \langle \psi_{p_j^{e_j}} \rangle = X_{2} \times
X^{(m)}$ where $m=n/2^{e}$.
\end{enumerate} 
\end{prop}

\begin{proof}
Part (a) is essentially part (a) of Lemma 1 in \cite{lettl}. Part (b) follows from the fact that 
the natural projection $X \rightarrow  \langle \psi_{2^{e}} \rangle$ and thus $X_{2} \rightarrow \langle \psi_{2^{e}} \rangle$ must be surjective.
By part (a), $\langle X, \prod_{i=2}^{t} \langle \omega_{p_{i}} \rangle \rangle$ contains all the Sylow-$p$ subgroups of $X^{(n)}= \prod_{i=1}^t \langle \omega_{p_i} \rangle \times \prod_{j=1}^t \langle \psi_{p_j^{e_j}} \rangle$ for $p$ odd; in particular, it contains 
$\prod_{j=2}^t \langle \psi_{p_j^{e_j}} \rangle$. Thus
$\prod_{i=2}^t \langle \omega_{p_i} \rangle \times \prod_{j=2}^t \langle \psi_{p_j^{e_j}} \rangle
 \subseteq \langle X , \prod_{i=2}^{t} \langle \omega_{p_{i}} \rangle \rangle$. Part (c) now follows by noting that the image of the natural projection $ \langle X, \prod_{i=2}^{t} \langle \omega_{p_{i}} \rangle \rangle \rightarrow
 X^{(2^{e})}$ is $X_{2}$. 
\end{proof}

\section{Ramification}

\begin{definition}
Throughout this paper, we take ``tamely ramified'' to mean ``at most tamely ramified'', i.e. ``not wildly ramified''.
\end{definition}

\begin{theorem}\label{wildtrace}
Let $L/K$ be an extension of number fields. Then $T_{L/K}(\mathcal{O}_L)$ is an ideal of $\mathcal{O}_K$. Suppose further that $L/K$ is Galois, and let $\mathfrak{p}$ be a (non-zero) prime of $\mathcal{O}_K$. Then $\mathfrak{p} \mid T_{L/K}(\mathcal{O}_L)$ if and only if $\mathfrak{p}$ is wildly ramified in $L/K$.
\end{theorem}

\begin{proof}
See \cite{mj}. Alternatively, this follows Lemma 2 in section 5 of \cite{lf} and the fact that the extension of residue fields in question is separable.
\end{proof}

\begin{corollary}\label{wildtracecor}
If $L/K$ is a Galois extension of number fields, then $L/K$ is tamely ramified
if and only if $T_{L/K}(\mathcal{O}_L)=\mathcal{O}_K$.
\end{corollary}

\begin{prop}\label{tameaboveodd}
Let $K$ be an abelian number field of conductor $n$. Then $\Q^{(n)}/K$ is tamely ramified at
each prime lying above an odd rational prime.
\end{prop}

\begin{proof}
Let $X$ be the group of Dirichlet characters associated to $K$ and write $n=\prod_{i=1}^t p_i^{e_i}$ where $p_{1}=2$. Let $M$ be the field corresponding to $\prod_{i=2}^{t} \langle \omega_{p_{i}} \rangle$. The extension
$MK/K$ is tamely ramified since the same is true of $M/ \Q$. By parts (b) and (c) of Proposition \ref{surjective} we have $[\Q^{(n)}:MK]=1$ or $2$, and so the result follows.
\end{proof}

\begin{corollary}\label{2subext}
Let $K$ be an abelian number field of conductor $n$. Then wild ramification in $\Q^{(n)}/K$ can only occur in a degree $2$ sub-extension (at primes above $2$). 
\end{corollary}

\begin{remark}
The result of Proposition \ref{tameaboveodd} appears to be well-known (it is noted
in \cite{byott_lettl}, for example), but its proof and corollary are not easily found in the literature.
\end{remark}

\begin{prop}\label{no2problems}
Let $K$ be an abelian number field of conductor $n=\prod_{i=1}^t p_i^{e_i}$ with associated 
character group $X$. Let $e=e_1=v_2(n)$. Then the following are equivalent:
\begin{enumerate}
\item $X_2 = X^{(2^{e})}$.
\item $X^{(n)} = \langle X, \prod_{i=2}^t \langle \omega_{p_i} \rangle \rangle$.
\item $\Q^{(n)} / K$ is tamely ramified.
\item $T_{\Q^{(n)}/K}(\mathcal{O}^{(n)}) = \mathcal{O}_{K}$, i.e. $I(\Q^{(n)}/K)=1$.
\end{enumerate}
\end{prop}

\begin{proof}
(a) $\Leftrightarrow$ (b) follows from part (c) of Proposition \ref{surjective}. \par
(c) $\Leftrightarrow$ (d) follows from Corollary \ref{wildtracecor}. \par
(a) $\Leftrightarrow$ (c) follows from Proposition \ref{tameaboveodd} and Theorem \ref{ramindex}.
\end{proof}

\begin{remark}
In particular, the equivalent conditions of Proposition \ref{no2problems}
hold when $e \leq 2$. Furthermore, it can be shown that they also hold if there exists $d \in \Z$ with $d \equiv 3 \, (4)$ and $d$ square-free such that $\Q[\sqrt{d}] \subseteq K$.
\end{remark}

\begin{prop}\label{some2problems}
Let $K$ be an abelian number field of conductor $n=\prod_{i=1}^t p_i^{e_i}$ with associated 
character group $X$ and let $K_{2}$ be the field corresponding to $X_{2}$. Let $e=e_1=v_2(n)$ and $m=n/2^{e}$. Define $L$ to be the compositum 
$K_{2}\Q^{(m)}$, i.e. the field corresponding to $X_{2} \times X^{(m)} \subseteq X^{(n)}$.
When the equivalent conditions of Proposition \ref{no2problems} do \emph{not} hold, the following statements are true:
\begin{enumerate}
\item $X_2$ is either $\langle \psi_{2^{e}} \rangle$ or $\langle \omega_{2} \psi_{2^{e}} \rangle$.
\item $L/K$ is tamely ramified.
\item $\Q^{(n)} = L[i]$, i.e. $\Q^{(n)}$ is the field generated by adjoining a root of $x^{2}+1$ to $L$.
\item $[\Q^{(n)}:L]=[L[i]:L]=2$.
\item $\Q^{(n)}/L$ is wildly ramified at the primes above 2.
\item $T_{L/K}(\mathcal{O}_{L}) = \mathcal{O}_{K}$.
\item $\mathcal{O}_{L}=\mathcal{O}_{K_{2}} \tensor_{\Z} \mathcal{O}^{(m)}$.
\item $I(\Q^{(n)}/L)=2^{r}$ for some $r \geq 1$.
\end{enumerate}
The situation is partially illustrated by the following field diagram.
$$
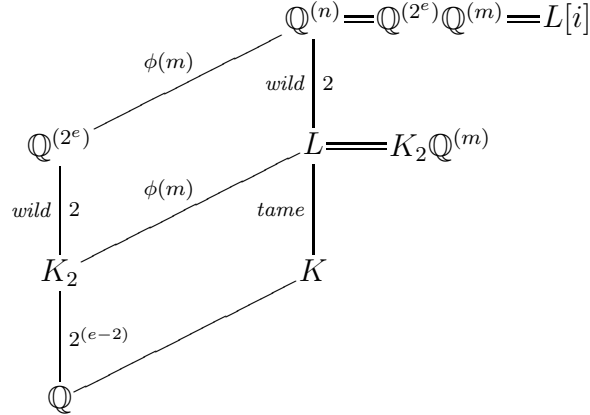
\xymatrix@1@!0@=48pt { 
&  & \Q^{(n)} \ar@{-}[d]_{\textrm{wild}}^{2} \ar@{=}[r] & \Q^{(2^{e})} \Q^{(m)} \ar@{=}[r] & L[i] \\ 
\Q^{(2^{e})} \ar@{-}[urr]^{\phi(m)} \ar@{-}[d]_{\textrm{wild}}^{2} & & L \ar@{-}[d]_{\textrm{tame}} 
\ar@{=}[r] & K_{2} \Q^{(m)}\\
K_{2} \ar@{-}[urr]^{\phi(m)} \ar@{-}[d]^{2^{(e-2)}} & & K \\
\Q \ar@{-}[urr]
}
$$
\end{prop}

\begin{proof}
(a) This follows from part (b) of Proposition \ref{surjective} and the hypothesis that part (a) of Proposition \ref{no2problems} does not hold.\par
(b) Since $X_{2} \times X^{(m)} = \langle X, \prod_{i=2}^{t} \langle \omega_{p_{i}} \rangle \rangle$ 
(Proposition \ref{surjective}, part (c)), the result follows by noting that $L=KM$ in the proof of 
Proposition \ref{tameaboveodd}.\par
(c) Since $\langle \omega_{2} \rangle$ corresponds to $\Q[i]$, this follows from part (a).\par
(d) $[\Q^{(n)}:L]=[X^{(n)}:X_{2} \times X^{(m)}] = [X^{(2^{e})}:X_{2}]=2$.\par
(e) This follows from part (b) and the hypothesis that part (c) of Proposition \ref{no2problems} does not hold (i.e. $\Q^{(n)}/K$ is wildly ramified). \par
(f) By Corollary \ref{wildtracecor}, this is equivalent to part (b). \par
(g) Since the discriminants of $\mathcal{O}_{K_{2}}$ and $\mathcal{O}^{(m)}$ are coprime, this follows from III.2.13 in \cite{ft}.\par
(h) This follows from part (e) and Theorem \ref{wildtrace}.
\end{proof}

\begin{prop}\label{ramcriteria}
Let $L/K$ be an extension of absolutely abelian number fields of equal conductor, $n$. Then each prime above an odd rational prime is tamely ramified in $L/K$. Furthermore, $L/K$ is wildly ramified at primes above $2$ if and only if:
\begin{enumerate}
\item the equivalent conditions of Proposition \ref{no2problems} applied to $L$ hold; and
\item the equivalent conditions of Proposition \ref{no2problems} applied to $K$ do \emph{not} hold.
\end{enumerate}
\end{prop}

\begin{proof}
Since $L/K$ is a sub-extension of $\Q^{(n)}/K$, the first statement follows from Proposition \ref{tameaboveodd}.
The second statement holds because wild ramification in $\Q^{(n)}/K$ can only occur in a degree $2$ sub-extension (Corollary \ref{2subext}), so $L/K$ is wildly ramified (at primes above 2) 
if and only if $\Q^{(n)}/L$ is tamely ramified and $\Q^{(n)}/K$ is wildly ramified.
\end{proof}

\section{Abelian number fields of conductor $2^{e}$, $e \geq 3$}

In this section, let $e \geq 3$, let $\zeta$ denote a primitive $2^{e}$-th root of unity and let $i = \zeta^{2^{e-2}}$.

\begin{prop}\label{ro2i}
The cyclotomic field $\Q^{(2^{e})}$ has precisely two proper fields of conductor $2^{e}$:
\begin{enumerate}
\item $\Q(\zeta+\zeta^{-1})$, with ring of integers $\Z[\zeta+\zeta^{-1}]$; and
\item $\Q(i(\zeta+\zeta^{-1}))$, with ring of integers $\Z[i(\zeta+\zeta^{-1})]$. 
\end{enumerate}
\end{prop}

\begin{proof}
Proposition \ref{surjective} part (b) implies that any proper subfield of $\Q^{(2^{e})}$ of conductor $2^{e}$ has associated character group either $\langle \psi_{2^{e}} \rangle$ or $\langle \omega_{2}  \psi_{2^{e}} \rangle$. It is straightforward to check that these correspond to $\Q(\zeta+\zeta^{-1})$ and $\Q(i(\zeta+\zeta^{-1}))$. \par
The ring of integers of $\Q(\zeta+\zeta^{-1})$ is $\Z[\zeta+\zeta^{-1}]$ by Proposition 2.16 of \cite{wash}. A slightly modified version of this argument, keeping track of real and imaginary parts, shows that $\Q(i(\zeta+\zeta^{-1}))$ has ring of integers $\Z[i(\zeta+\zeta^{-1})]$.
\end{proof}

\begin{prop}\label{image2trace}
Let $K_{2}$ be a proper subfield of $\Q^{(2^{e})}$ of conductor $2^{e}$. Let $T=T_{\Q^{(2^{e})}/K_{2}}$. 
In the cases of Proposition \ref{ro2i}, we have
\begin{enumerate}
\item $T(\Z[\zeta]) = 2\Z \oplus (\zeta+\zeta^{-1}) \cdot \mathcal{O}_{K_{2}} = 2\Z \oplus 
(\zeta+\zeta^{-1}) \cdot \Z[\zeta+\zeta^{-1}] $; and
\item $T(\Z[\zeta]) = 2\Z \oplus i(\zeta+\zeta^{-1}) \cdot \mathcal{O}_{K_{2}} = 2\Z \oplus 
i(\zeta+\zeta^{-1}) \cdot \Z[i(\zeta+\zeta^{-1})]$ .
\end{enumerate}
In both cases, $I(\Q^{(2^{e})}/K_{2})=2$.
\end{prop}

\begin{proof}
(a) In this case, $K_{2} = \Q(\zeta+\zeta^{-1})$, $\mathcal{O}_{K_{2}}= \Z[\zeta+\zeta^{-1}]$ and
$\{ 1, \zeta \}$ is a basis for $\Q^{(2^{e})}$ over $K_{2}$. The only non-trivial automorphism of 
$\Q^{(2^{e})}$ over $K_{2}$ is induced by complex conjugation, and so for $a,b \in K_{2}$, we have
$$T(a+ \zeta b) = (a+\zeta b) + (a +\zeta^{-1}b) = 2a + (\zeta+\zeta^{-1})b \, .$$
Since $\Z + \zeta \cdot \Z[\zeta+\zeta^{-1}] \subseteq \Z[\zeta]$, we therefore have
$ 2\Z \oplus (\zeta+\zeta^{-1}) \cdot \Z[\zeta+\zeta^{-1}] \subseteq T(\Z[\zeta])$.
However, $\Z[\zeta+\zeta^{-1}] = \Z \oplus (\zeta+\zeta^{-1}) \cdot \Z[\zeta+\zeta^{-1}]$, so
 $$\left[ \Z[\zeta+\zeta^{-1}]:2\Z \oplus (\zeta+\zeta^{-1}) \cdot \Z[\zeta+\zeta^{-1}] \right]=2$$ 
and by part (h) of Proposition \ref{some2problems}, 
$$\left[\Z[\zeta+\zeta^{-1}]:T(\Z[\zeta])\right]=2^{r}$$ for some $r \geq 1$.
Hence $2\Z \oplus (\zeta+\zeta^{-1}) \cdot \Z[\zeta+\zeta^{-1}] = T(\Z[\zeta])$ (and in fact $r=1$). \par
(b) In this case, $K_{2} = \Q(i(\zeta+\zeta^{-1}))$ and $\mathcal{O}_{K_{2}}= \Z[i(\zeta+\zeta^{-1})]$. The proof is essentially the same as in part (a), noting that $\{ 1 , i\zeta^{-1} = \zeta^{2^{(e-2)}-1}\}$ is a basis for $\Q^{(2^{e})}$ over $K_{2}$ and that the non-trivial Galois conjugate of $i\zeta^{-1}=\zeta^{2^{(e-2)}-1}$
over $K_{2}$ is $i\zeta=\zeta^{2^{(e-2)}+1}$.
\end{proof}

\begin{prop}\label{new2basis}
Consider the cases of Proposition \ref{ro2i}.
\begin{enumerate}
\item Let $A= \{ \zeta+\zeta^{-1}, \zeta^{2}+\zeta^{-2}, \ldots, \zeta^{2^{(e-2)}-1} + \zeta^{-2^{(e-2)}+1}\}$. \\ Then
$T(\Z[\zeta]) = \Span_{\Z}(A \cup \{2\})$, $\mathcal{O}_{K_{2}}=\Span_{\Z}(A \cup \{ 1 \})$
and $\Gal(K_{2}/\Q)(A) \subseteq \pm A$.
\item Let $B= \{ i(\zeta+\zeta^{-1}), \zeta^{2}+\zeta^{-2}, i(\zeta^{3}+\zeta^{-3}), \ldots,
i(\zeta^{2^{(e-2)}-1} + \zeta^{-2^{(e-2)}+1})\}$. \\ Then
$T(\Z[\zeta]) = \Span_{\Z}(B \cup \{2\})$, $\mathcal{O}_{K_{2}}=\Span_{\Z}(B \cup \{ 1 \})$ 
and $\Gal(K_{2}/\Q)(B) \subseteq \pm B$.
\end{enumerate}
\end{prop}

\begin{proof}
(a) $T(\Z[\zeta]) = \Span_{\Z}(A \cup \{2\})$ by Proposition \ref{image2trace} and a 
straight-forward induction argument; that $\mathcal{O}_{K_{2}}=\Span_{\Z}(A \cup \{ 1 \})$
follows easily. For any $\sigma \in \Gal(K_{2}/\Q)$ and any $j \in \{ 1, \ldots, 2^{(e-2)}-1 \} $, 
$\sigma(\zeta^{j}+\zeta^{-j}) = \zeta^{jk}+\zeta^{-jk}$
for some $k \in (\Z/2^{e}\Z)^{\times}$. However, any such $\zeta^{jk}+\zeta^{-jk}$ can be rewritten as
$\pm (\zeta^{r}+\zeta^{-r})$ for some $r \in \{ 1, \ldots, 2^{(e-2)}-1 \}$ (note $\zeta^{2^{(e-1)}}=-1$).
Part (b) is similar, noting that $\sigma(i) = \pm i$.
\end{proof}

\section{Computing $I(L/K)$}

\begin{prop}\label{ignoretame}
Let $L \subseteq M \subseteq N$ be a tower of Galois number fields such that $N/M$ is tamely ramified. Then $I(N/L) = I(M/L)$.
\end{prop}

\begin{proof}
Since $T_{N/L}(\mathcal{O}_{N}) = T_{M/L}(T_{N/M}(\mathcal{O}_{N}))$ and $T_{N/M}(\mathcal{O}_{N})=\mathcal{O}_{M}$ (by Corollary \ref{wildtracecor}), we have that
$T_{N/L}(\mathcal{O}_{N}) = T_{M/L}(\mathcal{O}_{M})$ and so the result follows from the definition of $I$.
\end{proof}

\begin{corollary}\label{ignoretamecor}
Let $L/K$ be a wildly ramified extension of absolutely abelian number fields of equal conductor, $n$.
Then $I(L/K) = I(\Q^{(n)}/K)$.
\end{corollary}

\begin{proof}
$\Q^{(n)}/L$ is tamely ramified since wild ramification in $\Q^{(n)}/K$ only occurs in a degree $2$ sub-extension (Corollary \ref{2subext}) and $L/K$ is wildly ramified.
\end{proof}

\begin{lemma}\label{tracetensor}
Let $K$ and $M$ be abelian number fields of conductors $k$ and $m$ respectively. Suppose that $k$ and $m$ are relatively prime. Then
$$ T_{\Q^{(k)}M/KM}(\mathcal{O}_{\Q^{(k)}M}) = T_{\Q^{(k)}/K}(\mathcal{O}^{(k)})
\tensor_{\Z} \mathcal{O}_{M} \, .$$
\end{lemma}

\begin{proof}
The proof is straightforward once one observes that by III.2.13 in \cite{ft}, we have
$\mathcal{O}_{KM}=\mathcal{O}_{K} \tensor_{\Z} \mathcal{O}_{M}$ and
$\mathcal{O}_{\Q^{(k)}M}=\mathcal{O}^{(k)} \tensor_{\Z} \mathcal{O}_{M}$.
\end{proof}

\begin{prop}\label{finalcompprop}
Let $K$ be an abelian number field of conductor $n$ such that $\Q^{(n)}/K$ is wildly ramified. Let $m=n/2^{e}$ where $e=v_{2}(n)$ and let $L=K_{2} \tensor_{\Q} \Q^{(m)} = K_{2}\Q^{(m)}$ (as in Proposition \ref{some2problems}). Let $C=A$ or $B$ from Proposition \ref{new2basis}, as appropriate. Define
$$ D = T_{L/K}(\mathcal{O}^{(m)}) \quad \textrm{and} \quad  E = T_{L/K}( \Span_{\Z}(C) \tensor_{\Z} \mathcal{O}^{(m)} ) \, . $$
Then 
$$\mathcal{O}_{K} = D \oplus E \quad \textrm{and} \quad T_{\Q^{(n)}/K}(\mathcal{O}^{(n)}) = 2D \oplus E \, .$$
\end{prop}

\begin{proof}
Note that $D \subseteq \mathcal{O}^{(m)} = \Z \tensor_{\Z} \mathcal{O}^{(m)}$ and
$E \subseteq \Span_{\Z}(C) \tensor_{\Z} \mathcal{O}^{(m)} $, with the last containment following from Proposition \ref{new2basis} (note $\Gal(L/K)(C) \subseteq \Gal(K_{2}/\Q)(C) \subseteq \pm C$). Since $\Z \cap \Span_{\Z}(C) = \{ 0 \}$, we have $D \cap E = \{ 0 \}$, which gives the last equality of 
\begin{eqnarray*}
\mathcal{O}_{K} &=& T_{L/K}(\mathcal{O}_{L}) \quad \textrm{(Proposition \ref{some2problems}, part (f))} \\
& = & T_{L/K}(\mathcal{O}_{K_{2}} \tensor_{\Z} \mathcal{O}^{(m)} )  
\quad \textrm{(Proposition \ref{some2problems}, part (g))} \\
&=& T_{L/K}((\Z \oplus \Span_{\Z}(C)) \tensor_{\Z} \mathcal{O}^{(m)} ) 
 \quad \textrm{(Proposition \ref{new2basis})} \\
& = & T_{L/K}((\Z \tensor_{\Z} \mathcal{O}^{(m)}) \oplus 
(\Span_{\Z}(C) \tensor_{\Z} \mathcal{O}^{(m)})) \\
&=& D+E = D \oplus E \, .
\end{eqnarray*}
Furthermore, 
\begin{eqnarray*}
T_{\Q^{(n)}/K}(\mathcal{O}^{(n)}) &=& T_{L/K}(T_{\Q^{(n)}/L}(\mathcal{O}^{(n)}))
= T_{L/K}(T_{\Q^{(n)}/L}(\mathcal{O}^{(2^{e})} \tensor_{\Z} \mathcal{O}^{(m)})) \\
&=& T_{L/K}(T_{\Q^{(2^{e})}/K_{2}}(\mathcal{O}^{(2^{e})}) \tensor_{\Z} \mathcal{O}^{(m)})
\quad \textrm{(Lemma \ref{tracetensor})} \\
&=& T_{L/K}((2\Z \oplus \Span_{\Z}(C)) \tensor_{\Z} \mathcal{O}^{(m)} )
\quad \textrm{(Proposition \ref{new2basis})} \\
&=& 2D \oplus E \quad \textrm{(as above)}.
\end{eqnarray*}
\end{proof}

\begin{remark}
The key point in this proof is the use of Proposition \ref{new2basis} to show that $D \cap E = \{ 0 \}$, and hence that the sums $D+E$ and $2D+E$ are direct.
\end{remark}

\begin{theorem}\label{finalcomp}
Under the hypotheses of Proposition \ref{finalcompprop}, we have 
\begin{enumerate}
\item $\mathcal{O}_{K} = \mathcal{O}_{K \cap \Q^{(m)}} 
\oplus T_{L/K}( \Span_{\Z}(C) \tensor_{\Z} \mathcal{O}^{(m)} )$;
\item $ T_{\Q^{(n)}/K}(\mathcal{O}^{(n)}) = 2\mathcal{O}_{K \cap \Q^{(m)}} 
\oplus T_{L/K}( \Span_{\Z}(C) \tensor_{\Z} \mathcal{O}^{(m)} )$; and
\item $ I(\Q^{(n)}/K) = 2^{[K \cap \Q^{(m)}:\Q]} .$
\end{enumerate}
\end{theorem}

\begin{proof}
Note that $\mathcal{O}_{K \cap \Q^{(m)}} \subseteq \mathcal{O}^{(m)} = \Z \tensor_{\Z} \mathcal{O}^{(m)}$ and, as shown in Proposition \ref{finalcompprop}, 
$E \subseteq \Span_{\Z}(C) \tensor_{\Z} \mathcal{O}^{(m)} $. Since $\Z \cap \Span_{\Z}(C) = \{ 0 \}$, we have $ \mathcal{O}_{K \cap \Q^{(m)}} \cap E = \{ 0 \}$
(this is essentially the same argument as that used to show $D \cap E = \{ 0 \}$). Furthermore,
$D = T_{L/K}(\mathcal{O}^{(m)}) \subseteq \mathcal{O}_{K \cap \Q^{(m)}}$ and 
$\mathcal{O}_{K \cap \Q^{(m)}} \subseteq \mathcal{O}_{K} = D \oplus E$, so $D = \mathcal{O}_{K \cap \Q^{(m)}}$. By Proposition \ref{finalcompprop}, this gives parts (a) and (b).
Now we have
\begin{eqnarray*}
 I(\Q^{(n)}/K) &=& [\mathcal{O}_{K} : T_{\Q^{(n)}/K}(\mathcal{O}^{(n)})] = 
[ \mathcal{O}_{K \cap \Q^{(m)}}\oplus E: 2\mathcal{O}_{K \cap \Q^{(m)}} \oplus E] \\
&=& [\mathcal{O}_{K \cap \Q^{(m)}}:2\mathcal{O}_{K \cap \Q^{(m)}}] = 2^{\rank_{\Z}(\mathcal{O}_{K \cap \Q^{(m)}})}
=2^{[K \cap \Q^{(m)}:\Q]},
\end{eqnarray*}
giving part (c).
\end{proof}

\begin{theorem}\label{mainresult}
Let $L/K$ be an extension of absolutely abelian number fields of equal conductor, $n$.
Let $e=v_{2}(n)$ and $m=n/2^{e}$. Then
$$
I(L/K) = \left\{ 
\begin{array}{ll} 
2^{[K \cap \Q^{(m)}:\Q]} =  2^{( [K:\Q] / 2^{(e-2)} )} & \textrm{ if $L/K$ is wildly ramified,} \\
1 & \textrm{ otherwise.}
\end{array} \right.
$$
\end{theorem}

\begin{remark}
Recall that criteria for wild ramification of $L/K$ (which can only happen at primes above $2$) are given in Proposition \ref{ramcriteria}.
\end{remark}

\begin{proof}
Suppose $L/K$ is wildly ramified. Then $I(L/K)=I(\Q^{(n)}/K)$ by Corollary \ref{ignoretamecor} and
$I(\Q^{(n)}/K) = 2^{[K \cap \Q^{(m)}:\Q]}$ by Theorem \ref{finalcomp}. Noting that $[K_{2}:\Q]=2^{(e-2)}$ 
(see Proposition \ref{some2problems}, part (a)) and that $\Q^{(m)}K=\Q^{(m)}K_{2}$, we have
\begin{eqnarray*}
[K \cap \Q^{(m)}:\Q] = \frac{[\Q^{(m)}:\Q]}{[\Q^{(m)}:K \cap \Q^{(m)}]}
= \frac{[\Q^{(m)}:\Q]}{[\Q^{(m)}K:K]}
= \frac{[\Q^{(m)}:\Q]}{[\Q^{(m)}K_{2}:K]} \\
= \frac{[\Q^{(m)}:\Q][K:\Q]}{[\Q^{(m)}K_{2}:\Q]}
= \frac{[\Q^{(m)}:\Q][K:\Q]}{[\Q^{(m)}:\Q][K_{2}:\Q]}
= \frac{[K:\Q]}{[K_{2}:\Q]}
= \frac{[K:\Q]}{2^{(e-2)}} \, .
\end{eqnarray*}
In the case where $L/K$ is tamely ramified, the result follows from Corollary \ref{wildtracecor}.
\end{proof}

\begin{remark}
It is clear that Theorem \ref{mainresult} agrees with the expressions for $I(L/K)$ in \cite{girstmair} (where $K \cap \Q^{(m)}$ is denoted $K_{n/2^{e}}$), and is in fact a sharpening of these results since an exact value for $I(L/K)$ is given for \emph{any} extension of abelian number fields $L/K$ of equal conductor. Furthermore, the above result does not rely on Leopoldt's Theorem.
\end{remark}

\section{The Adjusted Trace Map}

\begin{definition}\label{adjustdef}
Let $K$ be an abelian number field of conductor $n$. We define the ``adjusted trace map'', $\hat{T}_{\Q^{(n)}/K}$.
If $\Q^{(n)}/K$ is tamely ramified, let $\hat{T}_{\Q^{(n)}/K} = T_{\Q^{(n)}/K}$. Otherwise, let $m=n/2^{e}$ where $e=v_{2}(n)$ (recall that $e \geq 3$ in this case). Note that
$\mathcal{O}^{(n)}=\mathcal{O}^{(2^{e})} \tensor_{\Z} \mathcal{O}^{(m)}$ has $\Z$-basis 
$\{ \zeta_{2^{e}}^{i} \tensor \zeta_{m}^{j} \, \mid \, 0 \leq i \leq 2^{(e-1)}-1, \,  0 \leq j \leq \phi(m)-1 \}$. Define
$$
\hat{T}_{\Q^{(n)}/K}(\zeta_{2^{e}}^{i} \tensor \zeta_{m}^{j}) = \left\{ 
\begin{array}{ll} 
\frac{1}{2} T_{\Q^{(n)}/K}(\zeta_{m}^{j}) & \textrm{ for } i=0, \\
T_{\Q^{(n)}/K}(\zeta_{2^{e}}^{i} \tensor \zeta_{m}^{j}) & \textrm{ for } 1 \leq i \leq 2^{(e-1)}-1,
\end{array} \right. $$
and extend to a $\Q$-linear map $\Q^{(n)} \rightarrow K$.
\end{definition}

\begin{prop}\label{adjust_surjective}
$\hat{T}_{\Q^{(n)}/K}(\mathcal{O}^{(n)}) = \mathcal{O}_{K}$.
\end{prop}

\begin{proof}
If $\Q^{(n)}/K$ is tamely ramified, this is just Corollary \ref{wildtracecor}. Otherwise,
using the notation of Proposition \ref{finalcompprop}, we see that 
$$
\hat{T}_{\Q^{(n)}/K}(\alpha) = \left\{ 
\begin{array}{ll} 
\frac{1}{2} T_{\Q^{(n)}/K}(\alpha) & \textrm{ if } T_{\Q^{(n)}/K}(\alpha) \in D, \\
T_{\Q^{(n)}/K}(\alpha) & \textrm{ if }  T_{\Q^{(n)}/K}(\alpha) \in E.
\end{array} \right. 
$$
The result now follows immediately from Proposition \ref{finalcompprop}.
\end{proof}

\begin{remark}
It must be noted that the adjusted trace map of Definition \ref{adjustdef} is in fact equivalent to the
map defined in Lemma 3.4 of \cite{lux-pahlings} (page 51), though it is expressed more explicitly here. Furthermore, it is shown to be surjective in \cite{breuer}. However,
the proof given here (Proposition \ref{adjust_surjective}) is very different.
\end{remark}

\begin{lemma}\label{galinvariant}
Let
$\hat{T}_{\Q^{(n)}/K}(\zeta_{n}^{k}) = \varepsilon T_{\Q^{(n)}/K}(\zeta_{n}^{k})$ where $\varepsilon = 1/2$ or $1$. Then 
$$\hat{T}_{\Q^{(n)}/K}(\sigma(\zeta_{n}^{k})) = \varepsilon T_{\Q^{(n)}/K}(\sigma(\zeta_{n}^{k})) \quad \forall \sigma \in \Gal(\Q^{(n)}/\Q).$$
\end{lemma}

\begin{proof}
Write $\zeta_{n}^{k} = \zeta_{2^{e}}^{i} \tensor \zeta_{m}^{j}$ and use Definition \ref{adjustdef}.
\end{proof}

\begin{definition}
Let $L /K$ be a finite Galois extension with $G=\Gal(L / K)$.
Then $$\mathcal{A}_{L/K} := \{ \gamma \in K[G] \, | \, \gamma(\mathcal{O}_L) \subseteq \mathcal{O}_L \}$$
is the \emph{associated order} of $L / K$. 
\end{definition}

The following is a modified version of Lemma 6 in \cite{byott_lettl}. Note that we use both juxtaposition and the symbol
$\cdot$ to denote the action of a group algebra on a field.
 
\begin{theorem}
Let $K$ be an abelian number field of conductor $n$, and put $G = \Gal(\Q^{(n)}/\Q)$,
$H= \Gal(\Q^{(n)}/K)$. Let $\pi: \Q[G] \rightarrow \Q[G/H]$ denote the $\Q$-linear map induced by the natural projection $G \rightarrow G/H$.
Suppose $\mathcal{O}^{(n)}= \mathcal{A}_{\Q^{(n)}/\Q} \cdot \alpha$ for some $\alpha \in \mathcal{O}^{(n)}$. Then 
$\mathcal{A}_{K/\Q}= \pi(\mathcal{A}_{\Q^{(n)}/\Q})$ and $\mathcal{O}_{K}=\mathcal{A}_{K/\Q} \cdot \beta$
where $\beta = \hat{T}_{\Q^{(n)}/K}(\alpha)$.
\end{theorem}

\begin{proof}
Write $G= \{ g_{1}, \ldots, g_{r} \}$ and $H = \{ h_{1}, \ldots, h_{s}\}$.
Let $x \in \mathcal{A}_{\Q^{(n)}/\Q}$ and write
\begin{eqnarray*}
x &=& x_{1}g_{1} + \ldots + x_{r}g_{r} \textrm{ where } x_{i} \in \Q \textrm{ and }g_{i} \in G,\\
\alpha &=& y_{1} + y_{2}\zeta + \ldots + y_{r}\zeta^{r-1} \textrm{ where } y_{i} \in \Q \textrm{ and } \zeta=\zeta_{n}.
\end{eqnarray*}
Then using Lemma \ref{galinvariant}, the $\Q$-linearity of $\hat{T}_{\Q^{(n)}/K}$ and that $G$ is abelian, we have
\begin{eqnarray*}
\hat{T}_{\Q^{(n)}/K}(x \alpha) 
&=& \sum_{i=1}^{r} x_{i} \hat{T}_{\Q^{(n)}/K}(g_{i} \alpha )
= \sum_{i=1}^{r} x_{i} \sum_{j=1}^{r} y_{j} \hat{T}_{\Q^{(n)}/K}(g_{i} \zeta^{j-1}) \\
&=& \sum_{i=1}^{r} x_{i} \sum_{j=1}^{r} y_{j} \varepsilon_{j} T_{\Q^{(n)}/K}(g_{i} \zeta^{j-1})
= \sum_{i=1}^{r} x_{i} \sum_{j=1}^{r} y_{j} \varepsilon_{j} \sum_{k=1}^{s} h_{k}g_{i}\zeta^{j-1} \\
&=& \sum_{i=1}^{r} x_{i} g_{i} \sum_{j=1}^{r} y_{j} \varepsilon_{j} \sum_{k=1}^{s} h_{k}\zeta^{j-1}
= \sum_{i=1}^{r} x_{i} g_{i} \sum_{j=1}^{r} y_{j} \varepsilon_{j} T_{\Q^{(n)}/K}(\zeta^{j-1}) \\
&=& \sum_{i=1}^{r} x_{i} g_{i} \sum_{j=1}^{r} y_{j} \hat{T}_{\Q^{(n)}/K}(\zeta^{j-1})
= \sum_{i=1}^{r} x_{i} g_{i} \hat{T}_{\Q^{(n)}/K}(\alpha) \\ 
&=& x \hat{T}_{\Q^{(n)}/K}(\alpha) 
\end{eqnarray*}
where $\varepsilon_{j}=1/2$ or $1$, as appropriate. Thus
\begin{eqnarray*}
\mathcal{O}_{K} &=& \hat{T}_{\Q^{(n)}/K}(\mathcal{O}^{(n)})
\quad \textrm{(Proposition \ref{adjust_surjective})} \\
&=& \hat{T}_{\Q^{(n)}/K}( \mathcal{A}_{\Q^{(n)}/\Q} \cdot \alpha)
= \mathcal{A}_{\Q^{(n)}/\Q} \cdot \hat{T}_{\Q^{(n)}/K}(\alpha) \\
&=& \pi(\mathcal{A}_{\Q^{(n)}/\Q}) \cdot \beta 
\quad \textrm{ (since } \beta \in K \textrm{).}
\end{eqnarray*}
\end{proof}

\begin{remark}
Unfortunately, this result cannot be easily extended to the case of relative extensions because 
$\hat{T}_{\Q^{(n)}/K}$ is not $K$-linear for $K \neq \Q$.
\end{remark}

\begin{corollary}\label{reducetocyclotomic}
The proof of Leopoldt's Theorem can be reduced to the cyclotomic case.
\end{corollary}

We can now restate Leopoldt's Theorem (see \cite{leopoldt}, \cite{lettl}) with the generator expressed as the image of an element under the adjusted trace map.

\begin{definition}
For $n \in \N$, define the \emph{radical} of $n$ to be $r(n) = \prod_{p|n} p$. \end{definition}

\begin{definition}
For $n \in \N$, define $\mathcal{D}(n) = \{ d \in \N : r(n) | d  \textrm{ and } d | n \}$.
\end{definition}

\begin{theorem}[Leopoldt] 
Let $K$ be an abelian number field of conductor $n$, let $\zeta_{n}$ be a fixed primitive $n^{\mathrm{th}}$ root of unity, and let 
$$\alpha = \hat{T}_{\Q^{(n)}/K} \left( \sum_{d \in \mathcal{D}(n)} \zeta_{n}^{(n/d)} \right) \, .$$
Then we have $\mathcal{O}_K=\mathcal{A}_{K/\Q} \cdot \alpha$, and so $\mathcal{O}_K$ 
is a free $\mathcal{A}_K$-module of rank $1$. 
\end{theorem}

\begin{proof}
By Corollary \ref{reducetocyclotomic}, the proof is reduced to the cyclotomic case, which is relatively straightforward. 
\end{proof}

\begin{remark}
In particular, the cyclotomic case follows from the version of Leopoldt's Theorem given in \cite{lettl}.
\end{remark}

\begin{remark}
The definition of $\mathcal{D}(n)$ in \cite{lettl} is different from that given above. However, as noted in
\cite{lettl2}, Leopoldt's Theorem holds in either case. A routine computation shows that when $\mathcal{D}(n)$ is taken to be as in \cite{lettl}, $\alpha$ as defined above is equal to $T$ 
defined in \cite{lettl}.
\end{remark}  

\section{Acknowledgments}

The author is grateful to  Steven Chase, Ravi Ramakrishna, and Shankar Sen for useful conversations, and to Kurt Girstmair, Spencer Hamblen, and Jason Martin for looking at an initial draft of this paper. The computational algebra system Magma (\cite{magma}) was used to verify Theorem \ref{mainresult} for abelian number fields of conductor up to $176$. The positive results from this ``experiment'' were psychologically very helpful in proving the theorem.

\end{document}